\documentclass[a4paper]{amsart}
\usepackage{amssymb}
\usepackage{amscd}
\usepackage{amsthm}
\usepackage{amsmath}
\usepackage{mathtools}
\usepackage{latexsym}
\usepackage{hyperref}
\usepackage[all]{xy}
\usepackage[utf8]{inputenc}
\usepackage{enumitem}

\usepackage{tikz-cd}

\setlist[enumerate]{label={(\roman*)}}

\theoremstyle{plain}
\newtheorem{theorem}{Theorem}
\newtheorem{corollary}[theorem]{Corollary}
\newtheorem{lemma}[theorem]{Lemma}
\newtheorem{proposition}[theorem]{Proposition}

\theoremstyle{definition}

\newtheorem{example}[theorem]{Example}

\theoremstyle{remark}

\newtheorem{remark}{Remark}

\numberwithin{theorem}{section}
\usepackage{newtxtext}
\usepackage{newtxmath}

\DeclareMathAlphabet\urwscr{U}{urwchancal}{m}{n}%
\DeclareMathAlphabet\rsfscr{U}{rsfso}{m}{n}
\DeclareMathAlphabet\euscr{U}{eus}{m}{n}
\DeclareFontEncoding{LS2}{}{}
\DeclareFontSubstitution{LS2}{stix}{m}{n}
\DeclareMathAlphabet\stixcal{LS2}{stixcal}{m} {n}

\newcommand{\Spec}[1]{\operatorname{Spec}(#1)}

\newcommand{\card}{\mbox{\rm{card\,}}}

\newcommand{\add}{\mbox{\rm{add\,}}}

\newcommand{\Rej}[2]{\mbox{\rm{Rej}}_{#1}(#2)}
\newcommand{\Gen}{\mbox{\rm{Gen\,}}}
\newcommand{\Sum}{\mbox{\rm{sum\,}}}

\newcommand{\Cog}{\mbox{\rm{Cog\,}}}

\newcommand{\End}{\mbox{\rm{End\,}}}

\newcommand{\im}{\mbox{\rm{Im\,}}}

\newcommand{\Hom}[3]{\operatorname{Hom}_{#1}(#2,#3)}
\newcommand{\Ext}[4]{\operatorname{Ext}^{#1}_{#2}(#3,#4)}

\newcommand{\Tr}[2]{\mbox{\rm{Tr}}_{#1}(#2)}

\newcommand{\rmod}[1]{\mbox{\rm{Mod}--}{#1}}

\newcommand{\lmod}[1]{{#1}\mbox{--\rm{Mod}}}

\newcommand{\Ann}[1]{\mbox{\rm{Ann}}(#1)}

\begin{document}

\title{Dualizations of approximations, $\aleph_1$-projectivity, and Vop\v{e}nka's Principles}

\author{Asmae Ben Yassine and Jan Trlifaj}

\begin{abstract} The approximation classes of modules that arise as components of cotorsion pairs are tied up by Salce's duality. Here we consider general approximation classes of modules and investigate possibilities of dualization in dependence on closure properties of these classes. While some proofs are easily dualized, other dualizations require large cardinal principles, and some fail in ZFC, with counterexamples provided by classes of $\aleph_1$-projective modules over non-perfect rings. For example, we show that Vopěnka's Principle implies that each covering class of modules closed under homomorphic images is of the form $\Gen (M)$ for a module $M$, and that the latter property restricted to classes generated by $\aleph_1$-free abelian groups implies Weak Vop\v{e}nka's Principle.                
\end{abstract}

\date{\today}

\thanks{Research supported by GA\v CR 23-05148S. The first author also supported by SVV-2020-260589.}


\maketitle

\section*{Introduction}

Cotorsion pairs were introduced by Salce in \cite{Sal} as analogs of the well-known torsion pairs where the Hom functor was replaced by Ext. A formal replacement was certainly not the main point: Salce proved the remarkable fact that though there is no duality available in the category $\rmod R$, for each cotorsion pair $(\mathcal A,\mathcal B)$, the classes $\mathcal A$ and $\mathcal B$ are tied up by a duality: $\mathcal A$ is a special precovering class, if and only if $\mathcal B$ is a special preeveloping class, cf.\ \cite[Salce's Lemma 5.20]{GT}. 

For general classes $\mathscr C$ of modules that do not necessarily fit in cotorsion pairs, one still has the formally dual notions of a $\mathscr C$-preenvelope (or a left $\mathscr C$-approximation) and a $\mathscr C$-precover (or a right $\mathscr C$-approximation). However, there is no general tool for dualization like Salce's Lemma at hand. In the present paper, we consider general approximation classes of modules and investigate if, and how, dualizations are possible assuming extra closure properties of these classes.       

While some results can easily be dualized simply by employing dual arguments, other require the use of large cardinal principles. On the one hand, we prove that Vopěnka's Principle implies that each covering class of modules closed under homomorphic images is of the form $\Gen (M)$ for a module $M$ (Proposition \ref{p1}). On the other hand, we show that the latter property restricted to classes of abelian groups generated by $\aleph_1$-free groups implies Weak Vop\v{e}nka's Principle (Theorem \ref{t1}).  

In several cases, the class $\mathcal F \mathcal M$ of all $\aleph_1$-projective modules (= flat Mittag-Leffler modules) over a non-right perfect ring $R$ serves as a barrier for dualization in ZFC. While Weak Vop\v{e}nka's Principle is known to guarantee that each class of modules closed under direct products and direct summands is preenveloping (Lemma \ref{l1}), the dual statement is not true in ZFC: by Example \ref{e1}, $\mathcal F \mathcal M$ is a class of modules closed under direct sums and direct summands which is not precovering. Also, in contrast with the claim of Proposition \ref{p1} mentioned above, if $R$ is a Dedekind domain with a countable spectrum which is not a complete discrete valuation domain (e.g., when $R = \mathbb Z$), then $\mathcal F \mathcal M$ is an enveloping class of modules closed under submodules, but $\mathcal F \mathcal M$ is not of the form $\Cog (M)$ for any module $M$ (Example \ref{e2}).

\section{Preliminaries}

For a ring $ R $, we denote by $ \rmod R $ the class of all (right $ R $-) modules, and by $ \lmod R $ the class of all left $ R $-modules.
\subsection{Approximations}
A map $ f \in \Hom R {M} {C} $ with $ C \in \mathscr C  $ is a \emph{$ \mathscr C $-preenvelope} of $ M $, if the abelian group homomorphism $ \Hom R {f} {C'}:\Hom R {C}{C'} \rightarrow \Hom R {M}{C'}$ is surjective for each $ C' \in \mathscr C $.\\
	A $ \mathscr C $-preenvelope $ f \in \Hom R {M} {C}$ of $ M $ is called a \emph{$ \mathscr C $-envelope} of $ M $, provided that $ f $ is \emph{left minimal}, that is, provided $ f = gf $ implies that $ g $ is an automorphism for each $ g \in \End_{R}(\mathscr C) $.\\
	$ \mathscr C \subseteq \rmod R $ is a \emph{preenveloping class} (\emph{enveloping class}) provided that each module has a $ \mathscr C$-preenvelope ($\mathscr C $-envelope).

Dually, a map $ f \in \Hom R {C} {M} $ with $ C \in \mathscr C  $ is a \emph{$ \mathscr C $-precover} of $ M $, if the abelian group homomorphism $ \Hom R {C'} {f}:\Hom R {C'}{C} \rightarrow \Hom R {C'}{M}$ is surjective for each $ C' \in \mathscr C $.\\
	A $ \mathscr C $-precover $ f \in \Hom R {C} {M }$ of $ M $ is called a \emph{$ \mathscr C $-cover} of $ M $, provided that $ f $ is \emph{right minimal}, that is, provided $fg=f$ implies that $ g $ is an automorphism for each $ g \in \End_{R}(\mathscr C) $.\\
	$ \mathscr C \subseteq \rmod R $ is a \emph{precovering class} (\emph{covering class}) provided that each module has a $ \mathscr C$-precover ($\mathscr C $-cover). 

Let $ \mathscr C $ be a class of $ R $-modules. We define
$$ \mathscr C ^{\perp} = \{N \in \lmod{R} \;|\; \Ext{1}{R}{C}{N} = 0 \hbox{ for all } C \in \mathscr C\} $$
$$ ^{\perp}\mathscr C = \{N \in \lmod{R} \;|\; \Ext{1}{R}{N}{C}=0 \hbox{ for all } C \in \mathscr C\} $$

Let $ \mathcal A, \mathcal B \subseteq \rmod R $. The pair $ (\mathcal A, \mathcal B) $ is called a \emph{cotorsion pair} \cite[\S 5.2]{GT} (or a \emph{cotorsion theory}, \cite{Sal}), if $ \mathcal A = ^{\perp} \mathcal B $ and $ \mathcal B = \mathcal A^{\perp} $.

	A module $ M $ is called \emph{Enochs cotorsion} if $ \Ext{1}{R}{F}{M} = 0 $ for all flat modules $ F $. We denote by $ \mathcal{EC} $ the class of all Enochs cotorsion modules, and $\mathcal F _0$ the class of all flat modules. By \cite[Lemma 5.17]{GT}, $ (\mathcal F_{0}, \mathcal{EC}) $ forms a cotorsion pair, known as the \emph{Enochs cotorsion pair}.

A $\mathscr C$-preenvelope $f$ is called \emph{special} if $f$ is monic and its cokernel is an element of $^\perp \mathscr C$. A $\mathscr C$-precover $g$ is called \emph{special} if $g$ is surjective and its kernel is an element of $\mathscr C ^\perp$. A class $\mathscr C$ is \emph{special preenveloping} (\emph{special precovering}) in case each module has a special $\mathscr C$-preenvelope (special $\mathscr C$-precover).

A cotorsion pair $ (\mathcal A, \mathcal B) $ is called \emph{complete} if $\mathcal A$ is a special precovering class. By Salce's Lemma mentioned in the Introduction, this is equivalent to $\mathcal B$ being a special preenveloping class. 

\medskip
Throughout our paper, we stick to the terminology introduced for modules by Enochs \cite{E}. Notice that a different, though equivalent, terminology has been used in representation theory by the Auslander school (see e.g.\cite{AuR}): $\mathscr C$-preenvelopes and $\mathscr C$-envelopes are called \emph{left $\mathscr C$-approximations} and \emph{minimal left $\mathscr C$-approximations}. Dually, $\mathscr C$-precovers and $\mathscr C$-covers are called \emph{right $\mathscr C$-approximations} and \emph{minimal right $\mathscr C$-approximations}. 

In category theory, a still different terminology has been used, cf.\ \cite{AR}: Preenveloping classes closed under direct summands are called \emph{weakly reflective}, while  precovering classes closed under direct summands are \emph{weakly coreflective}. 

It is easy to see that all enveloping (covering) classes of modules are closed under direct summands, so they are  weakly reflective (weakly coreflective). An example of a precovering class that is not coreflective is provided by the class of all free modules over any ring possessing projective modules that are not free (such as a path algebra of a non-trivial quiver). Examples of preenveloping classes of modules that are not weakly reflective arise from pure-injective modules that are not dual modules: 

\begin{example}\label{ex1} Let $R$ be a ring, $\mathcal D$ be the class of all \emph{dual modules} (= the class of all modules isomorphic to the character modules $N^+ = \Hom{\mathbb Z}{N}{\mathbb Q/\mathbb Z}$ of all left $R$-modules $N$), and let $\mathcal P \mathcal I$ be the class of all \emph{pure-injective modules} (= modules injective w.r.t.\ all pure embeddings). By \cite[Theorem 2.27]{GT}, $\mathcal P \mathcal I$ is the class of all direct summands of the modules from $\mathcal D$. Both $\mathcal P \mathcal I$ and $\mathcal D$ are preenveloping classes of modules, since the canonical embedding of any module $M$ into its double dual module $M^{++}$ is pure, cf.\ \cite[Corollary 2.21(b)]{GT}. However, $\mathcal D \subsetneq \mathcal P \mathcal I$ in general:      

For an example, let $R = \mathbb Z$. Let $\mathcal E$ denote the class of all torsion-free divisible groups, that is, the underlying groups of $\mathbb Q$-linear spaces. By \cite[Corollary 2.18(a)]{GT}, ${\mathcal E}^+ = \mathcal D \cap \mathcal E$. However, $\mathbb Q \in \mathcal E \setminus {\mathcal E}^+$, so $\mathbb Q$ is a pure-injective abelian group which is not dual. The point is that duals of non-zero torsion-free divisible groups are uncountable, as $\mathbb Q^+$ is uncountable. Indeed, $\mathbb Q^+$ fits in the short exact sequence $0 \to (\mathbb Q/\mathbb Z)^+ \to \mathbb Q^+ \to \mathbb Z ^+ \to 0$ where $(\mathbb Q/\mathbb Z)^+ \cong \prod_{p \in \mathbb P} \mathbb J _p$ is an uncountable group. Here, $\mathbb P$ denotes the set of all prime integers and $\mathbb J _p$ the group of all $p$-adic integers, for $p \in \mathbb P$. 
\end{example}
\subsection{Modules}
	  Let $ R $ be a ring, $ M $ a module, and $ \mathscr C $ a class of modules. A family of submodules, $\mathcal{M}=(M_{\alpha} \mid \alpha \leq \sigma)$, of $ M $ is called a \textit{continuous chain} in $M$, provided that $M_{0}=0$, $M_{\alpha} \subseteq M_{\alpha+1}$ for each $\alpha < \sigma $, and $M_\alpha = \bigcup_{\beta < \alpha} M_{\beta}$ for each limit ordinal $\alpha \leq \sigma$.\\
	A continuous chain $\mathcal{M} $ in $M$ is a $ \mathscr C $-\textit{filtration} of $ M $, provided that $M=M_{\sigma}$, and each of the modules $M_{\alpha+1}/M_{\alpha}\;(\alpha < \sigma)$ is isomorphic to an element of $ \mathscr C $.\\
	$M$ is called $ \mathscr C $-\textit{filtered} (or a \emph{transfinite extension} of modules in $\mathscr C$), provided that $M$ possesses at least one $ \mathscr C $-filtration. A class $\mathscr C$ is \emph{closed under transfinite extensions} provided that $M \in \mathscr C$ for each $\mathscr C$-filtered module $M$. 
	
	Notice that each class of modules closed under transfinite extensions is closed under extensions and (arbitrary) direct sums.

\medskip
	Let $ \mathscr C $ be a class of modules. A module $ M $ is said to be \emph{generated} by $ \mathscr C $ if there exists a set $ \Lambda $, a family $ (C_{\lambda})_{\Lambda} $ of elements of $ \mathscr C $ and an epimorphism $\bigoplus_{\lambda \in \Lambda}C_{\lambda}\rightarrow M$.

Dually, let $ \mathscr C $ be a class of modules. A module $ M $ is said to be \emph{cogenerated} by $ \mathscr C $ if there exists a set $ \Lambda $, a family $ (C_{\lambda})_{\Lambda} $ of elements of $ \mathscr C $ and a monomorphism $M \rightarrow \prod_{\lambda \in \Lambda}C_{\lambda}$.

The class of all modules generated and cogenerated by $ \mathscr C $ is denoted by $ \Gen (\mathscr C) $ and $ \Cog (\mathscr C) $, respectively. If $ \mathscr C $ consists of a single module $ C $, we say that $ C $ generates (cogenerates) $ M $. We also use the notations $ \Gen (C)$ and $ \Cog (C) $ instead of $ \Gen (\mathscr C) $ and $ \Cog (\mathscr C )$.

We recall the following easy facts (see e.g.\ \cite[13.4 and 14.4]{Wis})

\begin{lemma}
	\begin{enumerate}
		\item The class $  \Gen(\mathscr C) $ is closed under homomorphic images and direct sums, and it contains $ \mathscr C $. Moreover if $ \mathcal X $ is a subclass of $ \rmod R $ which contains $ \mathscr C $, and is closed under epimorphic images and direct sums, then $ \Gen (\mathscr C) \subseteq \mathcal X$.
		\item The class $  \Cog(\mathscr C) $ is closed under submodules and direct products, and it contains $ \mathscr C $. Moreover if $ \mathcal X $ is a subclass of $ \rmod R $ which contains $ \mathscr C $, and is closed under submodules and direct products, then $ \Cog (\mathscr C) \subseteq \mathcal X$.
	\end{enumerate}
\end{lemma}

Let $M$ be a module. We will denote by $\Sum M$ the class of all finite direct sums of copies of $M$, and by $\add M$ the class of all direct summands of all modules in $\Sum M$. 

For a class of modules $\mathscr C$, we denote $\Tr{N}{\mathscr C}$ the \emph{trace} of $\mathscr C$ in $N$, that is, the sum of images of all homomorphisms from modules in $\mathscr C$ to $N$. Dually, $\Rej{N}{\mathscr C}$ denotes the \emph{reject} of $\mathscr C$ in $N$, that is, the intersection of kernels of all homomorphisms from $N$ to modules in $\mathscr C$, cf.\ \cite[p.109]{AF}. 

\begin{proposition}\label{pAF}\cite[Proposition 8.12]{AF}
	Let $ \mathscr C $ be a class of modules, and let $N$ be a module. Then:
	\begin{enumerate}
		\item $ \Tr{N}{\mathscr C} $ is the unique largest submodule $L$ of $N$ generated by $\mathscr C$;
		\item $\Rej{N}{\mathscr C}$ is the unique smallest submodule $U$ of $N$ such that $N/U$ is cogenerated by $\mathscr C$.
	\end{enumerate}
\end{proposition}

For a module $M$, $\sigma [M]$ denotes the class of all modules \emph{subgenerated} by the module $M$, that is, the submodules of all modules generated by $M$. This class is closed under submodules and homomorphic images, and it is the smallest Grothendieck subcategory of $\rmod R$ containing the module $M$, cf.\ \cite[\S 15]{Wis}. Dually, $\pi [M]$ will denote the class of all homomorphic images of modules cogenerated by $M$. Also this class is closed under homomorphic images and submodules. 

For a class of modules $\mathscr C$, we will denote by $\varinjlim \mathscr C$ the class of all modules that are direct limits of directed systems consisting of modules from $\mathscr C$. 

\subsection{$ \aleph_{1} $-projectivity}
	Let $ R $ be a ring and $ M $ be an $R$-module. We say that $ M $ is $ \aleph_{1} $-\emph{projective} in case there exists a set $ \mathcal S $ consisting of countably generated projective submodules of $ M $ with the following properties: $ 0 \in \mathcal S $, any countable subset of $ M $ is contained in an element of $ \mathcal S $, and $ \mathcal S $ is closed under unions of well-ordered chains of countable length.

Notice that if $R$ is a right hereditary ring, then a module is $ \aleph_{1} $-projective, iff each of its countably generated submodules is projective. In particular, $\aleph_1$-projective abelian groups (i.e., the abelian groups all of whose countable subgroups are free) are called $\aleph_1$-\emph{free}. 

	Let $ R $ be a ring. A module $ M $ is \textit{Mittag-Leffler} provided that the canonical group homomorphism
	$$ \varphi: M \otimes_R \prod_{i \in I} N_{i}\rightarrow \prod_{i \in I} M \otimes_R N_{i}$$
	defined by
	$$ \varphi(m \otimes_R (n_{i})_{i \in I})= (m \otimes_R n_{i})_{i \in I}$$
	is monic for each family $ (N_{i}\;|\; i \in I) $ of left $ R $-modules. \\
	Let $M \in \rmod R$ and $\mathcal Q \subseteq \lmod R$. Then $M$ is \emph{$\mathcal Q$-Mittag-Leffler}, provided that the canonical morphism $M \otimes {\prod_{i \in I} Q_i} \to \prod_{i \in I} M \otimes_R Q_i$ is injective for any family $( Q_i \mid i \in I )$ consisting of elements of $\mathcal Q$. So a module $M$ is Mittag-Leffler, iff it is $\mathcal Q$-Mittag-Leffler for $\mathcal Q = \lmod R$. 

We will be interested in flat Mittag-Leffler modules, and more in general, flat $\mathcal Q$-Mittag-Leffler modules. Following \cite{HT}, we will denote the class of all flat Mittag-Leffler modules by $\mathcal F \mathcal M$, and the class of all flat $\mathcal Q$-Mittag-Leffler modules by $\mathcal D _{\mathcal Q}$. 
 
The key relation between these notions goes back to \cite[Corollary 2.14(i)]{HT}: If $R$ is any ring and $M$ any module, then $M$ is $ \aleph_{1} $-projective, if and only if $M$ is flat Mittag-Leffler. In particular, abelian groups are flat Mittag-Leffler, iff they are $\aleph_1$-free.

\subsection{Vop\v{e}nka's Principles}
We will employ two large cardinal principles. The first one is due to Petr Vop\v{e}nka, cf.\ \cite[p.\ 278]{AR1}. It is now called \emph{Vop\v{e}nka's Principle}, and one of its equivalent renderings says that there exist no large rigid systems in the category $\mathcal G$ of all graphs. That is, there exists no proper class of graphs $\{ G_\alpha \mid \alpha \in \hbox{Ord} \}$ such that $\Hom{\mathcal G}{G_\alpha}{G_\beta} = \emptyset$ for all ordinals $\alpha \neq \beta$ and $\Hom{\mathcal G}{G_\alpha}{G_\alpha} = \{ \hbox{id}_{G_\alpha} \}$ for each ordinal $\alpha$. 

The second principle, called \emph{Weak Vop\v{e}nka's Principle}, says that there exists no proper class of graphs $( G_\alpha \mid \alpha \in \hbox{Ord} )$ such that for all ordinals $\alpha, \beta$, $\Hom{\mathcal G}{X_\alpha}{X_\beta} \neq \emptyset$, iff $\alpha \geq \beta$ \footnote{This is not the original formulation of Weak Vop\v{e}nka's Principle, but it is equivalent to it by \cite[Theorem 1.4]{W}}. 

It is known that Vop\v{e}nka's Principle implies Weak Vop\v{e}nka's Principle (\cite[Observation I.12]{AR}), but the converse fails under the assumption of existence of supercompact cardinals, see \cite[Theorem 1.2]{W}. A recent application of Vop\v{e}nka's Principle to approximation theory has appeared in \cite{C}: If Vop\v{e}nka's Principle holds, then each cotorsion pair over any right hereditary ring is complete. The latter claim cannot be proved in ZFC: by \cite{ES}, it is consistent with ZFC that the Whitehead cotorsion pair is not complete (the \emph{Whitehead cotorsion pair} is the cotorsion pair of abelian groups $(^\perp \mathbb Z,(^\perp \mathbb Z)^\perp)$). We refer to \cite[Appendix]{AR1} and \cite{W} for more details on the large cardinal strength of Vop\v{e}nka's principles.  

Also, we refer to \cite[Part II]{GT} for basics of approximation theory of modules, to \cite[Chap.\ IV and VII]{EM} for properties of $\aleph_1$-projective modules, and to \cite{AF} for basics of general theory of modules.  

\section{Closure properties, and enveloping and covering classes of modules}

We start with a characterization of preenveloping classes in terms of their closure properties:

\begin{lemma}\label{l1}
Let $R$ be a ring, and $\mathscr C \subseteq \rmod R$ a class of modules closed under direct summands. Consider the following two conditions:
\begin{enumerate}
		\item[(i)] $\mathscr C$ is preenveloping.
		\item[(ii)] $\mathscr C$ is closed under direct products.
\end{enumerate}
Then (i) implies (ii). If Weak Vop\v{e}nka's Principle holds, then (ii) implies (i). 
\end{lemma}
\begin{proof} Assume (i) and let $( E_i \mid i \in I )$ be a family of modules from $\mathscr C$. Let $f: P \to C$ be a $\mathscr C$-preenvelope of the module $P = \prod_{i \in I} E_i$. Denote by $\pi_i : P \to E_i$ the canonical projection ($i \in I$). Then there exist homomorphisms $g_i: C \to E_i$ such that $g_i f = \pi_i$ for each $i \in I$. Define a homomorphism $g : C \to P$ by $\pi_i g(c) = g_i(c)$ for all $c \in C$ and $i \in I$. Then $gf(x) = ( g_i(f(x)) \mid i \in I ) = x$ for all $x \in P$. Thus $P$ is isomorphic to a direct summand in $C$, and $P \in \mathscr C$ by our assumption on the class $\mathscr C$, so (ii) holds.

The implication (ii) implies (i) holds under Weak Vop\v{e}nka's Principle by \cite[Theorem 1.9 and Remark 1.10]{AR} and \cite[Theorem 1.4]{W}. 
\end{proof}

Next, we consider classes of modules closed under submodules. Under this additional assumption, the conditions (i) and (ii) of Lemma \ref{l1} become equivalent in ZFC:

\begin{lemma}\label{l1+}
Let $R$ be a ring, and $\mathscr C \subseteq \rmod R$ a class of modules. Consider the following conditions
\begin{enumerate}
		\item[(i)] $\mathscr C$ is (pre-) enveloping and closed under submodules.
		\item[(ii)] $\mathscr C$ is closed under submodules and direct products.
		\item[(iii)] $\mathscr C = \Cog (M)$ for a module $M$.
\end{enumerate}
Then (i) is equivalent to (ii), and it is implied by (iii). 
\end{lemma}
\begin{proof} (i) implies (ii) by Lemma \ref{l1}. Assume (ii) and let $ N \in \rmod R $. Let $U = \Rej{N}{\mathscr C}$. We claim that the canonical projection $ \pi_U: N \rightarrow N/U $ is a $ \mathscr C $-envelope of $ N $. 

First, $N/U \in \mathscr C$ by Proposition \ref{pAF}(ii). Let $ f \in \Hom RNC$ where $ C \in \mathscr C $. Since $\ker f \supseteq U$,   
the Homomorphism Theorem implies that $f$ factorizes through $ \pi_U $, and that the only factorization of $\pi_U$ through itself is by the identity automorphism $\hbox{id}_{N/U}$.    

Therefore, $\mathscr C$ is an enveloping class and (i) holds. Finally, (iii) trivially implies (ii).
\end{proof}

Let us try to dualize Lemma \ref{l1}. First, dualizing the proof of the implication (i) $\implies$ (ii) in Lemma \ref{l1}, we easily obtain

\begin{lemma}\label{l1d} 
Let $R$ be a ring, and $\mathscr C \subseteq \rmod R$ a class of modules closed under direct summands. Consider the following two conditions:
\begin{enumerate}
		\item[(i)] $\mathscr C$ is precovering.
		\item[(ii)] $\mathscr C$ is closed under direct sums.
\end{enumerate}
Then (i) implies (ii). 
\end{lemma}

The following two examples show that relative flat Mittag-Leffler modules yield a barrier in ZFC for proving both reverse implications, i.e., both (ii) $\implies$ (iii) in Lemma \ref{l1+}, and (ii) $\implies$ (i) in Lemma \ref{l1d}. We start with the second implication:

\begin{example}\label{e1} Let $R$ be a ring, $\mathcal Q$ a class of left $R$-modules, and $\mathcal D _{\mathcal Q}$ the class of all flat $\mathcal Q$-Mittag-Leffler modules So $\mathcal D _{\mathcal Q} = \mathcal F \mathcal M$ is the class of all flat Mittag-Leffler (= $\aleph_1$-projective) modules in the particular case when $\mathcal Q = \lmod R$, and 
$\mathcal D _{\mathcal Q} = \mathcal F \mathcal P$ the class of all $f$-projective modules in the particular case when $\mathcal Q = \{ R \}$, cf. \cite[\S 4]{BT}. By \cite[Lemma 3.1(i)]{BT}, the class $\mathcal D _{\mathcal Q}$ is always closed under transfinite extensions and pure submodules (and hence under direct sums and direct summands). 

By \cite[Theorem 3.6]{BT} (see also \cite[Theorem 3.3]{Sar}), for any class $\mathcal Q$ of left $R$-modules, the class $\mathcal D _{\mathcal Q}$ is precovering, iff $\mathcal D _{\mathcal Q}$ coincides with the class $\mathcal F _0$ of all flat modules. 

For $\mathcal Q = \lmod R$, it is well-known that $\mathcal F \mathcal M = \mathcal F_0$, only if the ring $R$ is right perfect, cf.\ \cite[\S 28]{AF}. So (ii) does not imply (i) in Lemma \ref{l1d} whenever $R$ is any non-right perfect ring and $\mathscr C = \mathcal F \mathcal M$. 

A different kind of examples arises for $\mathcal Q = \{ R \}$: by \cite[Proposition 4.8(ii)]{BT} and \cite[Proposition 4.10]{PR}, if $R$ is a right semihereditary ring such that $\mathcal F \mathcal P = \mathcal F _0$, then $R$ left semihereditary. So if $R$ is (the opposite ring of) the Chase ring from \cite[Chap.\ 1,\S2F,2.34]{L}, then (ii) does not imply (i) in Lemma \ref{l1d} for $\mathscr C = \mathcal F \mathcal P$.   
\end{example}

Before presenting our second example showing that in general (ii) does not imply (iii) in Lemma \ref{l1+}, we need a lemma generalizing a construction from \cite[Theorem 5.8]{EGPT}.

\begin{lemma}\label{noHom} Let $R$ be a right hereditary ring, $(\mathcal A, \mathcal B)$ be a cotorsion pair in $\rmod R$. Let $\mathscr C \subseteq \mathcal A$ be such that $\mathscr C$ is closed under submodules and transfinite extensions, $\mathscr C ^\perp = \mathcal B$, and $\mathscr C \cap \mathcal B = 0$.  

Then for each $0 \neq C \in \mathscr C$ there exists $0 \neq D \in \mathscr C$ such that $\Hom RDC = 0$. In particular, $\mathscr C \nsubseteq \Cog (C)$ for any module $C \in \mathscr C$.   
\end{lemma} 
\begin{proof} Let $0 \neq C \in \mathscr C$ and $\kappa = \card C + \card R + \aleph_0$. Let $\mathcal S$ denote a representative set of all non-zero modules in $\mathscr C$ of cardinality $\leq \kappa$ such that $C \in \mathcal S$. Since $\mathscr C ^\perp = \mathcal B$ and $\mathscr C \cap \mathcal B = 0$, for each $S \in \mathcal S$ there exists $C_S \in \mathscr C$ such that $\Ext 1R{C_S}S \neq 0$. Let $E = \bigoplus_{S \in \mathcal S} C_S \in \mathscr C$. Then $\Ext 1RES \neq 0$ for each $S \in \mathcal S$. Let $\lambda = \card{E} (\geq \kappa)$. Since $R$ is right hereditary, there is a projective resolution of $E$ of the form $0 \to K \overset{\eta}{\hookrightarrow} F \to E \to 0$ where $F$ is free of rank $\lambda$.

The module $D$ will be constructed as the last term of a $\mathscr C$-filtration $( D_\alpha \mid \alpha \leq \tau )$ for some $\tau \leq \lambda^+$ by induction as follows: $D_0 = C$; if $D_\alpha$ has already been constructed and $\Hom R{D_\alpha}C = 0$, we let $\tau = \alpha$ and $D = D_\tau$, and finish the construction.      

Otherwise we proceed by putting $H_\alpha = \Hom R{D_\alpha}C \setminus \{ 0 \}$. For each $h \in H_\alpha$, we let $I_h = \im h \subseteq C$. Since $\mathscr C$ is closed under submodules, $I_h \in \mathscr C$, and $\card I_h \leq \kappa$ implies $\Ext 1RE{I_h} \neq 0$ by our definition of $E$. Using the projective resolution of $E$ above, we infer that there exists $\phi_h \in \Hom RK{I_h}$ that cannot be extended to a homomorphism from $F$ to $I_h$. Since $K$ is projective and $h : D_\alpha \to I_h$ is surjective, there exists $\psi_h \in \Hom RK{D_\alpha}$ such that $\phi_h = h \psi_h$. 

We have the exact sequence $0 \to K^{(H_\alpha)} \overset{\eta_\alpha}{\hookrightarrow} F^{(H_\alpha)} \to E^{(H_\alpha)} \to 0$ where $E^{(H_\alpha)} \in \mathscr C$ because $\mathscr C$ is closed under transfinite extensions. For each $h \in H_\alpha$, let $\nu_h$ be the $h$th canonical inclusion of $K$ into $K^{(H_\alpha)}$, and $\mu_h$ $h$th canonical inclusion of $F$ into $F^{(H_\alpha)}$. Then $\mu_h \eta = \eta_\alpha \nu_h$ for each $h \in H_\alpha$. We define $\Psi_{\alpha} \in \Hom R{K^{(H_\alpha)}}{D_\alpha}$ by $\Psi_{\alpha} \nu_h = \psi_h$ for each $h \in H_\alpha$.  

The module $D_{\alpha + 1}$ is defined by the pushout of $\eta_\alpha$ and $\Psi_\alpha$,

$$\begin{CD}
0 @>>> K^{(H_\alpha)} @>{\eta_\alpha}>>  F^{(H_\alpha)} @>>> E^{(H_\alpha)} @>>> 0\\
@. @V{\Psi_\alpha}VV  @V{\rho}VV @| @.\\
0 @>>> D_\alpha @>{\subseteq}>> D_{\alpha + 1} @>>> E^{(H_\alpha)} @>>> 0.
\end{CD}$$

Then $D_{\alpha + 1} \in \mathscr C$ because $\mathscr C$ is closed under extensions. For a limit ordinal $\alpha$, we put $D_{\alpha} = \bigcup_{\beta < \alpha} D_{\beta}$. 

We claim that our construction terminates at some $\tau \leq \lambda^+$, whence for $D = D_\tau$, we have $\Hom RDC = 0$. 

If not, then $D_{\lambda^+}$ is defined and satisfies $D_{\lambda^+} = \bigcup_{\alpha < \lambda^+} D_{\alpha}$, and there exists $0 \neq g \in \Hom R{D_{\lambda^+}}C$. Let $\beta < \lambda^+$ be the least ordinal such that $g \restriction D_\beta \neq 0$. 

We will prove that for each $\beta \leq \alpha < \lambda^+$, $\im {(g \restriction D_\alpha)}$ is a proper submodule of $\im {(g \restriction D_{\alpha + 1})}$. Then $( \im {(g \restriction D_\alpha)} \mid \beta \leq \alpha < \lambda^+ )$ is a strictly increasing continuous chain of submodules of $C$ of cardinality $\lambda^+$, in contradiction with $\card C = \kappa \leq \lambda$.     

Assume $\im {(g \restriction D_\alpha)} = \im {(g \restriction D_{\alpha + 1})}$ for some $\beta \leq \alpha < \lambda^+$. Using the notation from the non-limit step of our construction for $h = g \restriction D_\alpha$, we define $f_h = (g \restriction D_{\alpha + 1}) \rho \mu_h \in \Hom RF{I_h}$. 

We claim that $f_h$ extends $\phi_h$ from $K$ to $F$, in contradiction with our definition of $\phi_h$. Indeed, $f_h \eta = (g \restriction D_{\alpha + 1}) \rho \eta_{\alpha} \nu_h$. Denoting the inclusion $D_\alpha \subseteq D_{\alpha + 1}$ by $\sigma_\alpha$ and using the pushout diagram above, we obtain  $f_h \eta = (g \restriction D_{\alpha + 1}) \sigma_\alpha \Psi_\alpha \nu_h$. The latter map is equal to $h \Psi_{\alpha} \nu_h$, and hence to $h \psi_h = \phi_h$. Thus, $f_h \eta = \phi_h$, q.e.d.
\end{proof}

\begin{example}\label{e2} Let $R$ be a Dedekind domain with a countable spectrum which is not a complete discrete valuation domain (e.g., let $R = \mathbb Z$). We claim that the class of all $\aleph_1$-projective modules $\mathcal F \mathcal M$ is an enveloping class closed under submodules, but $\mathcal F \mathcal M$ is not of the form $\Cog (M)$ for any module $M$.

To verify this claim, we apply Lemma \ref{noHom} to the Enochs cotorsion pair, i.e., $\mathcal A = \mathcal F_0$ and $\mathcal B = \mathcal{EC}$, and let $\mathscr C = \mathcal F \mathcal M$. This is possible, since $R$ is hereditary, whence $\mathscr C$ is closed under submodules by \cite[Corollary 2.10(i)]{HT}. Moreover $\mathscr C$ is closed under transfinite extensions by \cite[Corollary 3.20(i)]{GT}. Since $\Spec R$ is countable, the quotient field $Q$ is countably presented, whence $Q \in {}^\perp(\mathscr C^ \perp)$ by \cite[Lemma 10.19]{GT}. Then ${}^\perp(\mathscr C^ \perp)$ contains all torsion-free modules, whence ${}^\perp(\mathscr C^ \perp) = \mathcal F_0$, and $\mathscr C ^\perp = \mathcal B$.

Moreover, if $0 \neq M \in \mathscr C \cap \mathcal B$, then $M$ is a flat cotorsion module, so $M \cong Q^{(\alpha)} \oplus \prod_{0 \neq p \in \Spec{R}} \widehat{R_p ^{(\alpha_p)}}$ for some cardinals $\alpha$, $\alpha_p$ ($0 \neq p \in \Spec R$), where $R_p$ denotes the localization of $R$ at $p$ and \, $ \widehat{} $ \, the $p$-adic completion, see \cite[Theorem 5.3.28]{EJ}. Since $Q \notin \mathscr C$, $M$ has a direct summand isomorphic to $\widehat{R_p}$, whence $\widehat{R_p} \in \mathscr C$. We have an exact sequence $0 \to R_p \to \widehat{R_p} \to D \to 0$ where $D$ is torsion free and divisible, hence a direct sum of copies of $Q$. Since $R_p$ is a countably generated pure submodule of $\widehat{R_p}$, it is projective. This implies that $R = R_p$ is a (discrete) valuation domain. If $R$ is not complete, then there is an exact sequence $0 \to R \to N \to Q \to 0$ with $N \subseteq \widehat{R_p}$ countably generated, and hence free, a contradiction. Thus, if $R$ is not a complete discrete valuation domain, then $\mathscr C \cap \mathcal B = 0$.                    
  
Finally, since $R$ is noetherian, $\mathscr C$ is closed under direct products by \cite[Proposition 4.3]{HT}, so by Lemmas \ref{l1+} and \ref{noHom}, $\mathscr C$ is an enveloping class of modules closed under submodules, but $\mathscr C$ is not of the form $\Cog (M)$ for any module $M$.   
\end{example} 

We have seen that only one implication from Lemma \ref{l1} can be dualized. However, the dual version of Lemma \ref{l1+} does hold true, even in an extended form:

\begin{proposition}\label{p1}
Let $R$ be a ring, and $\mathscr C \subseteq \rmod R$ a class of modules. Consider the following conditions
\begin{enumerate}
		\item[(i)] $\mathscr C$ is (pre-) covering and closed under homomorphic images.
		\item[(ii)] $\mathscr C$ is closed under homomorphic images and direct sums.
		\item[(iii)] $\mathscr C = \Gen (M)$ for a module $M$.
\end{enumerate}
Then (i) is equivalent to (ii), and it is implied by (iii). If Vop\v{e}nka's Principle holds, then (ii) implies (iii).  
\end{proposition}
\begin{proof} (i) implies (ii) by Lemma \ref{l1d}. Assume (ii) and let $N \in \rmod R$. Let $T = T_{\mathscr C}(N)$. By (ii), $T \in \mathscr C$, and basic properties of the trace yield that the monomorphism $T \subseteq N$ is a $\mathscr C$-cover of $N$. So (i) holds. The implication (iii) $\implies$ (ii) is trivial.

Finally, assume (ii). Since $\mathscr C$ is closed under direct sums and homomorphic images, it is also closed under direct limits. Then Vop\v{e}nka's Principle gives $\mathscr C = \varinjlim \mathcal S$ for a subset $\mathcal S \subseteq \mathscr C$ (cf.\ \cite[Theorem 3.3]{ElB} or \cite[Corollary 6.18]{AR1}). Let $M = \bigoplus_{S \in \mathcal S} S$. Then $M \in \mathscr C$, and the closure properties of $\mathscr C$ yield that $\Gen (M)  \subseteq \mathscr C$. Moreover, by \cite[Lemma 1.1]{PPT}, $\Gen (M) \supseteq \varinjlim \Sum M = \varinjlim \mathcal S = \varinjlim \add M$, so $\Gen (M) = \mathscr C$, proving (iii).
\end{proof}

Adding further closure properties of the class $\mathscr C$ allows a complete characterization of the dual setting in ZFC in terms of the Grothendieck categories $\sigma[M]$. 

\begin{lemma}\label{l2}
Let $R$ be a ring, and $\mathscr C \subseteq \rmod R$ a class of modules. Then the following conditions are equivalent: 
\begin{enumerate}
		\item[(i)] $\mathscr C$ is (pre-) covering and closed under submodules and homomorphic images.
		\item[(ii)] $\mathscr C = \sigma [M]$ ($= \Gen (M)$) for a module $M$.
\end{enumerate}
\end{lemma}
\begin{proof} Assume (i). Denote by $\mathcal D$ a representative set of all finitely generated modules in the class $\mathscr C$. Let $M = \bigoplus_{D \in \mathcal D} D$. Since $\mathscr C$ is precovering and closed under direct summands, $\mathscr C$ is closed under direct sums by Lemma \ref{l1d}. So $M \in \mathscr C$. Then $\Gen (M) \subseteq \sigma [M] \subseteq \mathscr C$ as $\mathscr C$ is also closed under homomorphic images and submodules. Conversely, let $N \in \mathscr C$. Then $N$ is a directed union of copies of some modules from $\mathcal D$, say $N = \bigcup_{F \in \mathcal E} F$ for some $ \mathcal E \subseteq \mathcal D $. Thus $N \in \varinjlim \add M = \varinjlim \Sum M$ by \cite[Lemma 1.1]{PPT}, whence (ii) holds because $\varinjlim \Sum M \subseteq \Gen (M) \subseteq \sigma [M]$.

Assume (ii). Clearly, $\sigma [M]$ is closed under submodules and homomorphic images. For a module $N$, let $T = T_{\sigma[M]}(N)$. Then $T \in \sigma[M]$ and the inclusion $T \hookrightarrow N$ is a $\sigma [M]$-cover of $N$, whence (i) holds true.
\end{proof}

Let's turn again to the setting of preenveloping classes. For a module $M$, we will denote by $\pi [M]$ the class of all homomorphic images of all modules cogenerated by $M$.

\begin{lemma}\label{l2+}
Let $R$ be a ring, and $\mathscr C \subseteq \rmod R$ a class of modules. Then the following conditions are equivalent:
\begin{enumerate}
		\item[(i)] $\mathscr C$ is (pre-) enveloping and closed under submodules and homomorphic images.
		\item[(ii)] $\mathscr C = \pi [M]$ for a module $M$.
		\item[(iii)] $\mathscr C = \rmod {(R/I)}$ for a two-sided ideal $I$ in $R$.  
\end{enumerate}
\end{lemma}
\begin{proof} Assume (i). For each $M \in \mathscr C$, let $I_M = \Ann {M}$, and let $I = \bigcap_{M \in \mathscr C} I_M$. By Lemma \ref{l1}, $\mathscr C$ is closed under direct products, so $I = I_N$ for a module $N \in \mathscr C$, and moreover, there is a cyclic module $C = xR \in \mathscr C$ such that $I = \Ann {x}$. Thus $R/I \in \mathscr C$, and $\rmod {(R/I)} \subseteq \mathscr C$ by the closure properties of $\mathscr C$. However, $I \subseteq \Ann {M}$ for each module $M \in \mathscr C$, whence $\mathscr C \subseteq \rmod {(R/I)}$ and (iii) holds.       

Assume (iii). Then (ii) holds for $M = R/I \in \rmod R$. The proof that (ii) implies (i) is dual to the one given in Lemma \ref{l2}, replacing the trace by the reject. 
\end{proof}

\begin{remark}\label{r1} Since for each two-sided ideal $I$ in $R$, $\rmod {(R/I)} = \sigma [R/I]$, Lemmas \ref{l2} and \ref{l2+} imply that for each class $\mathscr C$ of modules closed under submodules and homomorphic images, $\mathscr C$ is covering whenever $\mathscr C$ is preenveloping. The converse fails in general as witnessed, for example, by the class $\mathscr C$ of all torsion abelian groups, cf.\ \cite[15.10]{W}. 

Notice also that for each module $M$, the class $ \sigma[M] $ is determined by the filter $ \mathcal{F}_{M}=\{I \leq R_{R} : R/I \in \sigma[M]\} $ consisting of right ideals of $ R $ (see e.g.\ \cite[\S 1]{PW}). The condition $ \sigma[M]= \rmod(R/I)$ for a two-sided ideal $I$ was characterized in \cite[Proposition 1.5]{PW}: $ \sigma[M]= \rmod(R/I)$, iff the filter $ \mathcal{F}_{M}$ is principal.     
\end{remark}

We finish by showing that the proof of the implication (ii) implies (iii) in Proposition \ref{p1} does in general require a large cardinal assumption, namely the Weak Vop\v{e}nka's Principle. The extra tool that we will need for this purpose is due to Prze\'{z}diecki, \cite[Theorem 3.14]{P}: 

\begin{lemma}\label{l3}    
There exists a functor $G$ from the category $\mathcal G$ of all graphs to $\rmod {\mathbb Z}$ which induces for all $X, Y \in \mathcal G$ a group isomorphism $\mathbb Z ^{({\Hom{\mathcal G}XY})} \cong \Hom{\mathbb Z}{G(X)}{G(Y)}$ natural in both variables.   
\end{lemma}

\begin{remark}\label{r2} By \cite[Corollary 4.11]{GP}, we can moreover assume that the functor $G$ from Lemma \ref{l3} takes its values in the class of all $\aleph_1$-free groups.
\end{remark}

\begin{theorem}\label{t1}
Assume that each covering class $\mathscr C$ of abelian groups which is closed under homomorphic images and is generated by a class of $\aleph_1$-free groups, satisfies $\mathscr C = \Gen (M)$ for an abelian group $M$. Then  Weak Vop\v{e}nka's Principle holds true. 
\end{theorem} 
\begin{proof} Assume that Weak Vop\v{e}nka's Principle fails, that is, there exists a proper class of graphs $( X_\alpha \mid \alpha \in \hbox{Ord} )$ such that for all ordinals $\alpha, \beta$, $\Hom{\mathcal G}{X_\alpha}{X_\beta} \neq \emptyset$, iff $\alpha \geq \beta$. 

Let $\mathscr C$ be the subclass of $\rmod {\mathbb Z}$ generated by the groups $G(X_\alpha)$ ($\alpha \in \mbox{Ord}$). By Remark \ref{r2}, we can w.l.o.g.\ assume that $G(X_\alpha)$ is $\aleph_1$-free for each $\alpha \in \mbox{Ord}$. Since $\mathscr C$ is closed under direct sums and homomorphic images, $\mathscr C$ is covering by Proposition \ref{p1}. We will show that there is no abelian group $M \in \mathscr C$ such that $\mathscr C = \Gen (M)$.

Assume that such a group $M$ does exist. Let $\alpha$ be the least ordinal such that $M$ is generated by the groups $G(X_\beta)$ ($\beta < \alpha$). Then $M$ is a homomorphic image of a direct sum of copies of those groups. Since $G(X_\alpha) \in \Gen (M)$, $G(X_\alpha)$ a homomorphic image of a direct sum of copies of $M$. Thus, there is a non-zero homomorphism from from $G(X_\beta)$ to $G(X_\alpha)$ for some $\beta < \alpha$. Then $\Hom{\mathcal G}{X_\beta}{X_\alpha} \neq \emptyset$ by Lemma \ref{l3}, a contradiction.      
\end{proof}

\begin{corollary}\label{c1} If Vop\v{e}nka's Principle holds, then each covering class of modules closed under homomorphic images is of the form $\Gen (M)$ for a module $M$. The latter property restricted to classes of abelian groups generated by $\aleph_1$-free groups implies Weak Vop\v{e}nka's Principle.                
\end{corollary}
\begin{proof} By Proposition \ref{p1} and Theorem \ref{t1}.
\end{proof}

\end{document}